\documentclass[11pt]{amsart}

\usepackage[T1]{fontenc}
\usepackage[utf8]{inputenc}
\usepackage{lmodern}

\usepackage[american]{babel}
\usepackage[babel, final]{microtype}
\usepackage{amssymb}
\usepackage[all]{xy}

\usepackage{ifthen}
\usepackage{ifpdf}
\usepackage{multirow}
\usepackage[shortlabels]{enumitem}
\usepackage{longtable}

\ifpdf
  \usepackage[pdftex, dvipsnames]{xcolor}
  \usepackage[pdftex, final]{graphicx}
  \usepackage{caption}
  \usepackage{subcaption}
  \usepackage{pinlabel}
  \usepackage[pdftex,%
    letterpaper,%
    includehead,%
    includefoot,%
    nomarginpar,%
    lmargin=1in,%
    rmargin=1in,%
    tmargin=1in,%
    bmargin=1in,%
  ]{geometry}
  \usepackage[pdftex,%
    final,%
    colorlinks=true,%
    linkcolor=NavyBlue,%
    citecolor=NavyBlue,%
    filecolor=NavyBlue,%
    menucolor=NavyBlue,%
    urlcolor=NavyBlue,%
    bookmarks=true,%
    bookmarksdepth=3,%
    bookmarksnumbered=true,%
    bookmarksopen=true,%
    bookmarksopenlevel=2,%
  ]{hyperref}
  \hypersetup{
    pdftitle={A Legendrian knot atlas for knots of arc index 10},
    pdfauthor={Ina Petkova, Noah Schwartz},
    pdfsubject={Legendrian knots},
    pdfkeywords={Legendrian knots}
  }
\else
  \usepackage[dvips, dvipsnames]{xcolor}
  \usepackage[dvips, final]{graphicx}
  \usepackage{caption}
  \usepackage{subcaption}
  \usepackage{pinlabel}
  \usepackage[dvips,%
    letterpaper,%
    includehead,%
    includefoot,%
    nomarginpar,%
    lmargin=1in,%
    rmargin=1in,%
    tmargin=1in,%
    bmargin=1in,%
  ]{geometry}
  \usepackage[ps2pdf,%
    final,%
    colorlinks=true,%
    linkcolor=NavyBlue,%
    citecolor=NavyBlue,%
    filecolor=NavyBlue,%
    menucolor=NavyBlue,%
    urlcolor=NavyBlue,%
    bookmarks=true,%
    bookmarksdepth=3,%
    bookmarksnumbered=true,%
    bookmarksopen=true,%
    bookmarksopenlevel=2,%
  ]{hyperref}
  \hypersetup{
    pdftitle={A Legendrian knot atlas for knots of arc index 10},
    pdfauthor={Ina Petkova, Noah Schwartz},
    pdfsubject={Legendrian knots},
    pdfkeywords={Legendrian knots}
  }
\fi

\usepackage{mathdots}
\usepackage{float}
\usepackage[mathscr]{euscript}
\usepackage[normalem]{ulem}
\usepackage{tikz}
\usepackage{tikz-cd}
\usepackage{mathtools}
\usepackage{stmaryrd}
\SetSymbolFont{stmry}{bold}{U}{stmry}{m}{n}
\usetikzlibrary{arrows,automata,fit}
\usetikzlibrary{matrix,arrows}
\usetikzlibrary{decorations.pathreplacing,calligraphy}

\input{macros}

\newcommand{\numset}[1]{\mathbb{#1}}

\newcommand{\R}{\numset{R}}

\newcommand{\xistd}{\xi_{\mathrm{std}}}
\DeclareMathOperator{\tb}{tb}
\DeclareMathOperator{\rot}{r}
\DeclareMathOperator{\self}{sl}

\newcommand{\curves}[1]{\boldsymbol{#1}}
\newcommand{\alphas}[1][]{%
  \ifthenelse{\equal{#1}{}}{\curves{\alpha}}{\curves{\alpha^{#1}}}
}
\newcommand{\betas}[1][]{%
  \ifthenelse{\equal{#1}{}}{\curves{\beta}}{\curves{\beta_{#1}}}
}

\newcommand{\grid}{\mathbb{G}}

\newcommand{\markers}[1]{\mathbb{#1}}
\newcommand{\OO}{\markers{O}}

\newcommand*{\alphastd}{\alpha_{\mathrm{std}}}

\renewcommand*{\tb}{\mathit{tb}}
\renewcommand*{\rot}{\mathit{r}}

\newcommand{\X}{\mathbb{X}}
\newcommand{\G}{\mathbb{G}}

\newcommand{\Leg}{\Lambda}

\newcommand{\std}{\mathrm{std}}

\newcommand*{\cref}{\fullfref}

\title[A Legendrian knot atlas for knots of arc index 10]{A Legendrian knot atlas for knots of arc index 10}

\author[I. Petkova]{Ina Petkova}
\address{Department of Mathematics \\ Dartmouth College \\ Hanover, NH 03755}
\email{\href{mailto:ina.petkova@dartmouth.edu}{ina.petkova@dartmouth.edu}}
\urladdr{\url{https://math.dartmouth.edu/~ina/}}

\author[N. Schwartz]{Noah Schwartz}
\address{Johns Hopkins University Applied Physics Laboratory \\ Laurel, MD 
  20723}
\email{\href{mailto:Noah.Schwartz@jhuapl.edu}{Noah.Schwartz@jhuapl.edu}}

\begin{document}

\begin{abstract}
  We expand the atlas of Legendrian knots  in standard contact three-space to knots of arc index 10. 
\end{abstract} 

\maketitle

\section{Introduction}\label{sec:intro}

One of the major problems in contact knot theory is that of classifying Legendrian and transverse knots of a fixed topological type $K$ in a given contact three-manifold $(Y, \xi)$. Even the simplest case where $Y = \R^3$ and $\xi = \xistd$ is the standard contact structure given by the kernel of the $1$-form
$  \alphastd = dz - y \,dx$
is quite challenging, and few Legendrian and transverse knots are fully classified; fully classified knots include the unknot \cite{EF-08}, torus knots \cite{EH01}, and twist knots \cite{ENV13}, as well as some satellite knots \cite{T-12, ELT-12}. 

In  2010, Chongchitmate and Ng undertook this task for all knots of arc index up to 9 \cite{CN13atlas}. The resulting ``Legendrian knot atlas'' contains conjectural classification data for these knots, and is updated and permanently available at {\url{https://services.math.duke.edu/~ng/atlas/}}.

In this paper, we begin the conjectural Legendrian and transverse classification  for prime knots of arc index 10. 
As in \cite{CN13atlas}, we take a probabilistic approach to enumerating Legendrian knots, by presenting them as grid diagrams. Two grid diagrams represent the same Legendrian knot if and only if they can be related by a sequence of elementary moves, which we call ``Legendrian grid moves''. One of these moves, stabilization, increases the size of the grid diagram. Thus, an algorithm cannot affirm in finite time that two grid diagrams represent nonisotopic Legendrian knots. Still, it is reasonable to guess that if two diagrams are related by Legendrian grid moves, then the sequence of moves does not increase grid size by much. In \cite{CN13atlas}, the conjectural classification is based on generating all grid diagrams of size 9 and all grid diagrams of size 10 (for knots of arc index 9). Here, we generate all (up to Legendrian isotopy) grid diagrams of size 10 for knots of (smooth) arc index 10.

 Starting with a large list of grid diagrams, we used a variety of methods to reduce it, by discarding duplicates -- at any given point, if two grids in the set represent Legendrian isotopic knots, we may remove one of them, and we still have a set that contains all Legendrian representatives that the original set contained. We guess that among the remaining diagrams that we present here, most pairs for a given topological knot are not Legendrian isotopic. Some pairs can indeed be distinguished by nonclassical invariants such as the ruling invariant or linearized contact homology, or the GRID invariants; we include an example in \fullref{sec:results}.
 We present the resulting data in the form of tables, as the graphical representation of grid diagrams and mountain ranges becomes somewhat unwieldly as the number and complexity of knots increases.

The atlas for knots of arc index 10 is available at \url{https://github.com/ipetkova/LegendrianAtlas}. Future updates to the data will be posted there. For example, we have also generated a set of grid diagrams of size 11 for these knots, and are running reduction code to eliminate duplicates (up to Legendrian isotopy). It is a lengthy process, but we hope to update the atlas with this data soon. The code we used to generate and reduce grids is also available at the above link.  In addition to that code, we used a combination of SnapPy and Sage methods to identify knots \cite{snapPy, sagemath}. 

While this manuscript was in preparation, the authors learned that a research group at the Louisiana State University had started work on a similar project. It would be interesting to compare our results with theirs when they become available \cite{lsu}.

\subsection*{Acknowledgments}

This work began in the Summer 2022, as part of an effort to generate interesting examples of filtered GRID invariants for our 2022 Summer Hybrid Undergraduate Research (SHUR) project. We thank Mike Wong for helpful conversations, and  Mitchell Jubeir for his contributions to efficiently generating an initial set of grid diagrams of size 10. We also thank Chuck Livingston for his help in filling gaps in identifying the topological types of some of our knots.Thanks to the anonymous referees for many helpful suggestions; in particular, thanks to the referee who brought to our attention the work of Dynnikov and Prasolov \cite{DP13} and its relevance to \cite{CN13atlas}. Much of the data in this paper appears in the second author's undergraduate honors thesis at Dartmouth College.
IP was partially supported by NSF CAREER Grant DMS-2145090. 

\section{Preliminaries}\label{sec:prelim}

In this section, we review some basics about Legendrian and transverse knots, as well as grid diagrams.

\subsection{Legendrian and transverse knots}

In this section, we review some basics of Legendrian  and transverse knots. We restrict to the case of knots in  $\R^3$ with  the standard contact structure 
\[
  \xi_{\std}=\ker(\alpha_{\std}),\quad \alpha_{\std}=dz-y\,dx.
\]
A smooth link $\Leg\subset(\R^3, \xi_{\std})$ is called \emph{Legendrian} if it is everywhere tangent to the contact structure.
Two Legendrian links are Legendrian isotopic if they are isotopic through a family of Legendrian links.

\begin{figure}[ht]
      \includegraphics[scale=1]{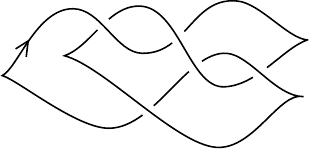}
  \caption{An example of a front projection.}
    \label{fig:4_1}
\end{figure}

A Legendrian link can be represented by its \emph{front diagram}, that is, by its projection onto the $xz$-plane. Note that in a front diagram, strands with lower slope always pass over strands with higher slope. See \fullref{fig:4_1} for an example. 

\begin{figure}[ht]
      \includegraphics[scale=1]{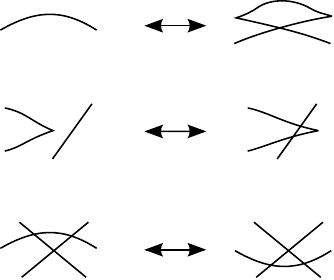}
  \caption{Legendrian Reidemeister moves. The vertical and horizontal reflections of the second move are also allowed.}
    \label{fig:leg-moves}
\end{figure}

Two Legendrian links are Legendrian isotopic if and only if their front diagrams are related by a sequence of  Legendrian planar isotopies (isotopies that preserve left and right cusps) and Legendrian Reidemeister moves (see
\fullref{fig:leg-moves}).

The two classical Legendrian link invariants are the Thurston--Bennequin number 
$\tb(\Leg)$ and the rotation number $\rot(\Leg)$. These can be computed from an 
oriented front diagram $D$ via the relations
\[
  \tb(\Leg)=\mathrm{wr}(D) - \frac{1}{2} (c_+(D)+c_-(D)), \qquad 
  \rot(\Leg)=\frac{1}{2} (c_-(D)-c_+(D)),
\]
where $\mathrm{wr}(D)$ is the writhe of the diagram, and $c_-(D)$ and $c_+(D)$ 
are the number of downward and upward cusps, respectively.

A smooth link  $ T \subset(\R^3, \xi_{\std})$ is called \emph{transverse} if it is everywhere transverse to the contact structure. Two transverse links are transversely isotopic if they are isotopic through transverse links. A transverse link $T$ naturally inherits an orientation from the (oriented) contact structure: a vector $v$ tangent to $T$ is positive if and only if $\alpha_{\std}(v)>0$.

Given an oriented Legendrian link $\Leg$, we can obtain a transverse link $T_+(\Leg)$, called the \emph{positive transverse pushoff of $\Leg$}, by smoothly pushing $\Leg$ off in a direction transverse to the contact planes so that the orientation of the smooth link is preserved. Further, every transverse link is the positive transverse pushoff of some Legendrian link. 

The \emph{self-linking number $\self(T)$} of a transverse link  $ T \subset(\R^3, \xi_{\std})$  is a classical transverse invariant, which can be defined as follows. If $\Leg$ is a Legendrian link such that $T$ is the positive transverse pushoff of $\Leg$, then 
\[\self(T) = \tb(\Leg) - r(\Leg).\]

Transverse links can be studied via front projections similar to those for Legendrian links. Instead of discussing those projections too, we shift the exposition to grid diagrams, the method of choice in this paper of representing both Legendrian and transverse links. 

\subsection{Grid diagrams}
\label{ssec:prelim-grid}

In this section, we review some basics of grid diagrams, following the 
conventions in \cite{OSS15} and \cite{OST08}.

A \emph{grid diagram} $\G$ is an $n\times n$ grid on the plane, along with two sets of markings
\[\OO = \{O_1, \ldots, O_n\}, \qquad \X = \{X_1, \ldots, X_n\},\]
such that there is exactly one $O$ and exactly one $X$ in each row, as well as in each column, and no square of the grid contains more than one marking. The number $n$ is called the \emph{grid number} of $\G$.

A grid diagram $\G$ specifies a link $L(\G)$ in $\R^3$ as follows. In each column, connect the $X$-marking to the $O$-marking with an oriented vertical segment. In each row, connected the $O$-marking  to the $X$-marking with an oriented horizontal segment.  Change resulting intersections to crossings, so  that vertical segments always cross above horizontal ones. We say \emph{$\G$ is a grid diagram for $L$}. Conversely, every link $L$ in $\R^3$ can be represented by a grid diagram.
By a theorem of Cromwell \cite{Cro95}, two grid diagrams represent the same link if and only if they are related by a sequence of the following grid moves: 
\begin{itemize}
\item \emph{cyclic permutations}, in which the topmost row (or bottommost row, rightmost column, leftmost column) is moved to become the bottommost row (or topmost row, leftmost column, rightmost column, resp.);
\item  \emph{commutations}, in which two adjacent rows (or columns) are switched if the corresponding segments connecting the $X$'s and $O$'s are either nested or disjoint;
\item \emph{stabilizations}, in which a $1 \times 1$ square with an $O$-marking (resp.\ $X$-marking)  is replaced by a $2 \times 2$ square with two diagonal $O$-markings and an $X$-marking (resp.\ two diagonal $X$-markings and an $O$-marking), creating a new row and a new column, and \emph{destabilizations}, the inverse operations to stabilizations. 
\end{itemize}
Following \cite{OSS15}, we classify (de)stabilizations by the type of the marking and the location of the empty cell in the $2 \times 2$ square; for example, a stabilization of type \textit{X:SE} results in a $2 \times 2$ square with an empty southeastern cell, an $X$ in the northwestern cell, and $O$'s in the other two cells.

A grid diagram $\G$ also specifies a Legendrian link $\Leg (\G)$ in $(\R^3, \xi_{\std})$, as follows. First, create the oriented link $L(\G)$.  The 
projection of $L(\G)$ onto the grid has corners that can be classified into 
four types: northeast, northwest, southwest, and southeast. Smooth 
all of the northwest and southeast corners of the projection, and turn the 
northeast and southwest corners into cusps. Next, rotate the diagram $45$ 
degrees clockwise. Last, flip all the crossings.
Note that the smooth type of the Legendrian link $\Leg (\G)$ is the mirror of the smooth link $L(\G)$. Similar to the smooth case, every Legendrian link in $(\R^3, \xistd)$ can be represented by a grid diagram. Two grid diagrams represent the same Legendrian link if and only if they are related by a sequence of cyclic permutations, commutations, and (de)stabilizations of type \textit{X:SE} and \textit{X:NW}. (De)stabilizations of type \textit{O:SE} and \textit{O:NW} also result in Legendrian isotopies. 
The other grid stabilization moves do change the Legendrian isotopy class of the link: 
\textit{X:NE} and \textit{O:SW} result in \emph{positive Legendrian stabilizations} and change $(\tb, r)$ by $(-1, 1)$, whereas \textit{X:SW} and \textit{O:NE} result in \emph{negative Legendrian stabilizations} and change $(\tb, r)$ by $(-1, -1)$.

If a Legendrian link is the positive/negative Legendrian stabilization of another, it is called \emph{destabilizable}; otherwise, it is called \emph{non-destabilizable}.
For a given topological type, $\tb$ is bounded above by a classical result of Bennequin \cite{B-83}, so there are non-destabilizable Legendrian links in every topological type.

A grid diagram $\G$ also specifies a transverse link $T(\G)$ in $(\R^3, \xi_{\std})$, by taking the positive pushoff of $\Leg(\G)$. Two grid diagrams represent the same transverse link if and only if they are related by a sequence of cyclic permutations, commutations, and (de)stabilizations of type \textit{X:SE}, \textit{X:NW}, and \textit{X:SW}. (De)stabilizations of type \textit{O:SE},  \textit{O:NW}, and \textit{O:SW} result in transverse isotopies as well. 
Thus, one can think of transverse links up to transverse isotopy as Legendrian links up to Legendrian isotopy and negative Legendrian stabilization. Stabilizations of type \textit{X:NE} and \textit{O:SW} result in \emph{transverse stabilization} of the transverse link $T(\G)$, and decrease the self-linking number by two. 

\begin{remark}
Our convention for converting from grid diagrams to Legendrian and transverse knots is different from \cite{CN13atlas}, and agrees with \cite{OSS15}, \cite{OST08}, and \cite{NOT08}.
\end{remark}

\subsection{Mountain ranges}

In  \cite{CN13atlas}, Legendrian representatives of a given topological knot type are depicted via a \emph{Legendrian mountain range}, as follows. Isotopy classes of Legendrian representatives are plotted as dots on the $(r, \tb)$-plane. Positive and negative Legendrian stabilizations are depicted by arrows. Black dots represent Legendrian representatives that are known to be distinct, while red dots represent conjecturally distinct representatives.

\subsection{Knot symmetries}

As in  \cite{CN13atlas}, we consider certain Legendrian symmetries. 

The \emph{reverse} $-\Leg$ of a Legendrian link  $\Leg$ is the result of  reversing the orientation of $\Leg$. On a grid diagram, this can be achieved by exchanging the $X$'s and  the $O$'s. This operation preserves $\tb$ and negates $r$. This operation may also change the topological type of the knot;  all knots in  \cite{CN13atlas} are isotopic to their orientation reverses, so this does not happen there.

The \emph{Legendrian mirror} $\mu(\Leg)$ of $\Leg$ is the image of $\Leg$ under the contactomorphism $(x,y,z)\mapsto (-x,y,-z)$. If $\G$ is a grid diagram for $\Leg$, one can obtain a grid diagram for $\mu(\Leg)$ by rotating $\G$ by $180^{\circ}$. This operation too  preserves $\tb$ and negates $r$; in addition, it preserves topological type.

The composition of the two operations above descends to an operation on transverse links, called the \emph{transverse mirror}.

\section{Methodology}

Our methodology is similar to  \cite{CN13atlas}. We briefly recall the main idea and outline major differences below.

The set of all grid diagrams can be thought of as an infinite graph $\Gamma$. Vertices correspond to grid diagrams; if two vertices are related by a grid move that corresponds to Legendrian isotopy, they are connected by an edge. The connected components of $\Gamma$ are precisely the isotopy classes of Legendrian links.\footnote{Note that if two vertices are in the same connected component, then they have the same topological type.} Denote by $\Gamma_n$ the finite subgraph corresponding to grid diagrams of size at most $n$. Determining the connected components of $\Gamma_n$ is the same as determining which pairs of vertices are connected by a path, and approximates the problem of finding the connected components of the infinite graph $\Gamma$.

We generate a set of grid diagrams of size $10$ that is guaranteed to contain a vertex in each connected component of $\Gamma_{10}$. While in \cite{CN13atlas} grid diagrams were represented by matrices of zeros and ones, we find it more efficient to represent grids by pairs of permutations, corresponding to the placement of the $X$'s and $O$'s. We discard grid diagrams representing multicomponent links, and subdivide the remaining grids into ``buckets'' by topological type, $\tb$, and $r$; to identify the topological type of a knot,  we used a combination of SnapPy and Sage methods\cite{snapPy, sagemath}, and consulted Chuck Livingston for help to fill the few remaining gaps \cite{chuck}. For knots that are not invertible (i.e. knots where reversing orientation changes the knot), we did not further subdivide by orientation, as identifying the knot with orientation can be a bit tricky. However, we were able to generate partial data, and used it to aid reduction.  For each bucket, we run a search to check if a grid can be destabilized, and if so, we discard it. We then run a bidirectional search to determine which pairs of grids are connected in $\Gamma_n$, and discard duplicates;  we started with $n=10$ (that is, we first searched for pairs of grids that are connected via a sequence of grid moves that does not involve grid stabilizations), and, as the set decreased in size, we check up to $n=12$ (that is, we searched for pairs of grids in $\Gamma_{10}$ that are connected in $\Gamma_{12}\supset\Gamma_{10}$). The goal is to reduce to a set of grids that conjecturally represents pairwise nonisotopic Legendrian knots. 
 
To manage time and space challenges, we sometimes relied on random sampling, imposed time limits on the code, and varied the depth of the bidirectional search and the maximum allowed size of the neighborhoods it generates. We also added edges to $\Gamma_{10}$ corresponding to additional moves that preserve Legendrian isotopy type and grid size, such as the $S_2$ move from \cite{NT08}.

\begin{remark}
We have also generated a set of  (conjecturally Legendrianly  nondestabilizable) grid diagrams of size $11$ for knots of (smooth) arc index $10$. We are running bidirectional searches to eliminate duplicates and reduce the size of this set, and will update the atlas one the data size is sufficiently small. 
\end{remark}

\section{Observations and results}\label{sec:results}

The authors of \cite{CN13atlas} conjectured that a Legendrian knot of maximal Thurston--Bennequin number has a grid diagram of minimum grid number, that is, a grid diagram of size the arc index of the knot. More generally, they conjectured that, given a knot $K$ of maximal Thurston--Bennequin number $\overline{\tb}(K)$ and arc index $\alpha(K)$, a Legendrian representative $\Leg$ for $K$ with $\tb(\Leg) = \overline{\tb}(K) - m$ can be represented on a grid of size $\alpha(K) + m$. As kindly pointed out by the referee, subsequent work of Dynnikov and Prasolov \cite{DP13} implies that this is indeed the case:

\begin{proposition}\label{prop:CN13-conj2}
Let $\Leg$ be a Legendrian representative for a knot $K$ with $\tb(\Leg) = \overline{\tb}(K) - m$.  Then $\Leg$ can be represented by a grid of size $\alpha(K)+m$, and this is the minimal size of a grid that represents $\Leg$. 
\end{proposition}
\begin{proof}
For any grid $\grid$, let $n(\grid)$ denote its size. 

In \cite{M06}, Matsuda shows that $\alpha(K)\geq \overline{\tb}(K) -  \overline{\tb}(m(K))$. We lay out a slightly generalized proof in the next paragraph. 

Let $\Leg$ be a Legendrian representative for a knot $K$ with $\tb(\Leg) = \overline{\tb}(K) - m$.  Suppose $\grid$ is a grid  for $\Leg$, and let $\grid'$ be the grid obtained by rotating $\grid$ clockwise by $90^{\circ}$. Observe that the corners of $\grid$ that become cusps of $\Leg(\grid)$ are exactly the complement of the corners of $\grid$ that become cusps of $\Leg(\grid')$, and that the positive/negative crossings of $\grid$ become negative/positive crossings in $\grid'$, so 
\[n(\grid) = -\tb(\Leg(\grid)) - \tb(\Leg(\grid'))\geq  - \overline{\tb}(K) -  \overline{\tb}(m(K))+ m.\]
 It follows that any grid for $\Leg$ must have size at least  $- \overline{\tb}(K) -  \overline{\tb}(m(K))+ m$. In particular, $\alpha(K)\geq - \overline{\tb}(K) -  \overline{\tb}(m(K))$.

Next, we show that this minimal grid size can always be realized. 

Let $\Leg$ be a Legendrian representative for a knot $K$ with $\tb(\Leg) = \overline{\tb}(K) - m$, and choose a Legendrian representative $\Leg'$ for $m(K)$ with $\tb(\Leg') = \overline{\tb}(m(K))$. By \cite[Theorem 7]{DP13}, there exists a grid diagram $\grid$ that realizes $\Leg$ and whose  $90^{\circ}$ clockwise rotation $\grid'$ realizes $\Leg'$. Again using that $n(\grid) = -\tb(\Leg(\grid)) - \tb(\Leg(\grid'))$, we get that 
\[n(\grid) = - \overline{\tb}(K) -  \overline{\tb}(m(K))+ m.\]

Note that in particular we get $\alpha(K) = - \overline{\tb}(K) -  \overline{\tb}(m(K))$. This completes the proof of the proposition.
\end{proof}

 \fullref{prop:CN13-conj2} immediately upgrades \cite[Conjecture 8]{CN13atlas} to a theorem, thus proving the existence of Legendrian mountain ranges with non-maximal peaks.
\begin{theorem}
There exist Legendrian mountain ranges with non-maximal peaks. In particular, for each of the knots $m(10_{139})$, $10_{161}$, and $m(12n242)$, there exists a Legendrian representative $\Leg$ with $\tb(\Leg)$ strictly less that the maximal possible $\tb$, and for which there is no other Legendrian representative $\Leg'$ with $(\tb(\Leg'), r(\Leg')) = (\tb(\Leg) + 1, r(\Leg)+1)$ or $(\tb(\Leg'), r(\Leg')) = (\tb(\Leg) + 1, r(\Leg)-1)$.
\end{theorem}
\begin{proof}
This is immediate from \fullref{prop:CN13-conj2}, since the list of Legendrian representatives realizable on grids of size 9 is exhaustive in \cite{CN13atlas}. 
\end{proof}

Further, \ref{prop:CN13-conj2} implies that for knots of arc index 10, all Legendrian representatives with maximal $\tb$ are listed in \fullref{tab:layer1}:

\begin{theorem}
There are at most $2,686$ Legendrian representatives with maximal $\tb$ and topological arc index $10$.  
\end{theorem}
\begin{proof}
This too is immediate from \fullref{prop:CN13-conj2}, since for knots of arc index 10 our list of Legendrian representatives realizable on grids of size 10 is exhaustive. 
\end{proof}

With orientation and mirroring, there are $604$ topological knots of arc index $10$ ($240$ knots, of which $65$ are chiral, $16$ negative amphicheiral, $161$ reversible, and $6$ fully amphicheiral). Of these $604$ knots, the program guesses that $135$ are transversely nonsimple:
 
 \begin{equation*}
\begin{matrix}
	m8a4, 8a6, m8a6, m8a7, 8a8, m8a8, 8a9, 8a10, m8a10, 8a11, 8a12, 8a13, \\ 
    8a14, m8a14,  8a15, m8a15, 8a16, m8a16, \pm 8a17,  10n4, m10n4, 10n6, 10n7, 10n8, m10n8, \\ 
    10n10, m10n10, m10n11, 10n12, m10n12, m10n14, 10n15, m10n15, m10n20, \\
    10n21, m10n23, \pm m10n24, \pm 10n25, \pm m10n25,  \pm 10n27, \pm m10n27, \pm 10n28,  \\
    m10n31, 10n32, 10n33, m10n33, m10n35, m10n36,  m10n39,   m10n40,\\
    m10n41, 10n42,  11n12, 11n20, \pm m11n23, \pm 11n24, \pm m11n24,  \\
    \pm m11n37, m11n48, \pm m11n50, 11n57, 11n61, m11n61, m11n65, 11n70, m11n82, \\
    11n84, m11n84, 11n86,   11n92, m11n96, m11n99, 11n106, m11n106, \\
    \pm 11n110, m11n111, m11n117,  \pm m11n122, \pm 11n132, \pm m11n132, \pm m11n134,  11n145, m11n164, \\
    12n25, 12n121, m12n243,  12n293, 12n321, 12n323, 12n340, m12n340, \\
    \pm 12n358, 12n370, m12n370, m12n375, 12n403,  m12n407, m12n438, \\
    12n443, 12n451, m12n452, 12n487, m12n487, m12n502, m12n603, 12n830, \\
     \pm 13n588, \pm 13n1692, 13n1907, \pm m13n1945, \pm m13n2787, \\
    13n3582, 15n41131.
\end{matrix}
\end{equation*}

Among the rest are the twist and torus knots of arc index $10$, namely $K_6 = 8a1$, $K_{-6} = m8a1$, $T_{3,7} = 14n21881$, and $T_{3,-7} = m14n21881$, all four of which are known to be transversely simple \cite{EH01, ENV13}.

\begin{figure}[ht]
      \includegraphics[scale=0.37]{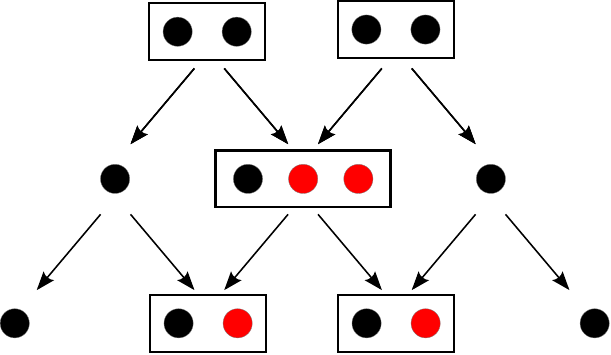}
  \caption{The Legendrian mountain range for the knot $12n121$. Black dots represent Legendrian representatives that are known to be distinct, while red dots represent conjecturally distinct representatives. 
}
    \label{fig:12n121}
\end{figure}

One could use non-classical invariants to distinguish some of the pairs of grids we list in this article. Examples of such invariants include the GRID invariants, their recently introduced  spectral generalizations \cite{JPSWW22}, the ruling invariant or linearized contact homology. For all but seven topological knot types, we have computed the GRID invariants of all their Legendrian representatives of maximal $\tb$, and we have provided those computations on \url{https://github.com/ipetkova/LegendrianAtlas}. For these fully-computed cases, the GRID invariants were only able to distinguish Legendrian representatives for the knot $12n121$:
\begin{proposition}
In the knot type $12n121$, there are exactly two distinct Legendrian representatives with $(\tb, r) = (0,-1)$ and two with $(\tb, r) = (0,1)$.
\end{proposition}
\begin{proof}
That there are at most two Legendrian representatives with $(\tb, r) = (0,-1)$ and at most two with $(\tb, r) = (0,1)$ follows from our exhaustive generation process, along with \ref{prop:CN13-conj2}. That all these representatives are distinct can be shown using the GRID invariants $\lambda^{\pm}$. Let $\Leg_1$, $\Leg_2$, $\Leg_3$, $\Leg_4$ be the Legendrian knots represented by grids G523, G524, G525, and G526, respectively.  A computation shows that $\lambda^-(\Leg_1)= 0$ and $\lambda^-(\Leg_1)\neq 0$, so $\Leg_1$ and $\Leg_2$ are not Legendrian isotopic. Similarly, $\lambda^+(\Leg_3)\neq 0$ and $\lambda^+(\Leg_4) = 0$, so $\Leg_3$ and $\Leg_4$ are not Legendrian isotopic. 
\end{proof}
Figure \fullref{fig:12n121} depicts the  mountain range for $12n121$. 

Computing the spectral GRID invariants has been slower, and we have fewer completed computations;  In the cases where we were able to compute the spectral GRID invariants, they were not able to distinguish additional pairs of grids.

\section{The Legendrian knot atlas for knots of arc index $10$}

The tables in this section present our conjectural classification of Legendrian knots of arc index $10$ (without yet considering nonmaximal Legendrian nondestabilizable knots). As in \cite{CN13atlas}, we present the first three horizontal layers of a conjectural mountain range for each prime knot $K$, that is, the Legendrian knots with Thurston--Bennequin numbers $\overline{\tb}(K)$, $\overline{\tb}(K)-1$, and $\overline{\tb}(K)-2$. Each of these layers, along with symmetry data, is presented in a corresponding table, as follows. 

\fullref{tab:layer1} contains a complete list of Legendrian nondestabilizable knots of arc index $10$ of maximal $\tb$. The list consists of grid diagram representatives. To track data more easily, we label each grid diagram with a unique identifier (its ID), of form ``G'' followed by a number. 
For each knot, we also list known symmetry data, as follows. In the row for grid $\G$, we list a grid $\G'$ in the column labeled ``reverse'', ``Legendrian mirror'', or ``reverse of Legendrian mirror'', if we can verify that $\Leg(\G')$ is Legendrian isotopic to $-\Leg(\G)$, $\mu(\Leg(\G))$, or $-\mu(\Leg(\G))$, respectively.  The goal of course is to determine whether a knot $\Leg$ is Legendrian isotopic to  $-\Leg$, $\mu(\Leg)$, and $-\mu(\Leg)$, but we find it helpful to include data rather than a question mark, even when $\Leg(\G')$ may be distinct from $\Leg(\G)$. For a  knot $\Leg$ with nonzero rotation number, we lists a hyphen in the columns for $-\Leg$ and $\mu(\Leg)$, as $-\Leg$ and $\mu(\Leg)$ are trivially distinct from $\Leg$. In all other cells of the table, if no $\G'$ as above is found, we mark the cell with a question mark.

\fullref{tab:layer2} contains a list of the positive and negative Legendrian stabilizations of the maximal knots. The grid IDs for these grids begin with an ``S''. 
For each grid, a list of ``parents'' is also given; these are the grids from \fullref{tab:layer1} whose Legendrian stabilization results in the given grid. Whether the respective Legendrian stabilization is positive or negative can be determined from the rotation number, so we omit listing this information for brevity. Note that by symmetry, we know that two of the four positive Legendrian stabilizations for $m8a14$ in  \fullref{tab:layer2} must be Legendrian isotopic, but we haven't been able to determine which ones yet. 

\fullref{tab:layer3}  contains a list of the  twice-Legendrian-stabilized maximal knots. Their grid IDs begin with a ``T'', and symmetry and ancestry data is listed in the same was as in \fullref{tab:layer2}. Note that for many of the knots, the fact that $S_+(S_-(\Leg)) = S_-(S_+(\Leg))$, along with the data from \fullref{tab:layer1}  and \fullref{tab:layer2}, completely determines the data in \fullref{tab:layer3}. We include the complete table nevertheless.

\clearpage
\footnotesize

\setlength{\tabcolsep}{4pt}

\begin{longtable}{|c|c|c|c|r|r|c|c|c|}
    \hline
 \textbf{Knot} & \textbf{ID} & \textbf{$\X$-permutation} & \textbf{$\OO$-permutation} & {$\boldsymbol{\tb}$} &  $\boldsymbol{r}$ & $\boldsymbol{\mu(\Leg)}$ & $\boldsymbol{-\Leg}$ & $\boldsymbol{-\mu(\Leg)}$ \\ \hline
\endhead
 \hline
        8a1 & G1011 & (0, 9, 2, 1, 4, 3, 7, 8, 5, 6) & (8, 1, 0, 3, 2, 5, 4, 6, 7, 9) & -7 & 0 & G1011 & G1011 & G1011 \\ \hline
m8a1 & G1848 & (0, 8, 1, 9, 4, 2, 5, 3, 7, 6) & (5, 2, 7, 3, 0, 6, 1, 8, 9, 4) & -3 & 0 & G1848 & G1848 & G1848 \\ \cline{2-9}
~ & G1849 & (0, 1, 8, 9, 6, 7, 5, 4, 2, 3) & (2, 9, 0, 7, 8, 4, 3, 6, 5, 1) & -3 & 0 & G1849 & G1849 & G1849 \\ \cline{2-9}
~ & G1850 & (0, 6, 4, 8, 7, 2, 1, 3, 5, 9) & (8, 9, 7, 6, 3, 5, 4, 0, 1, 2) & -3 & 0 & G1850 & G1850 & G1850 \\ \hline
8a2 & G1014 & (0, 2, 3, 9, 8, 7, 6, 5, 4, 1) & (3, 4, 1, 2, 0, 9, 8, 7, 6, 5) & -11 & -4 & -- & -- & G1014 \\ \cline{2-9}
~ & G1015 & (0, 9, 8, 7, 4, 5, 3, 2, 6, 1) & (2, 1, 0, 9, 6, 8, 7, 4, 3, 5) & -11 & -4 & -- & -- & G1015 \\ \cline{2-9}
~ & G1012 & (0, 9, 7, 8, 6, 5, 3, 2, 4, 1) & (2, 1, 0, 5, 4, 9, 7, 6, 8, 3) & -11 & -2 & -- & -- & G1013 \\ \cline{2-9}
~ & G1013 & (0, 7, 6, 4, 5, 3, 2, 1, 9, 8) & (6, 5, 9, 8, 0, 7, 4, 3, 2, 1) & -11 & -2 & -- & -- & G1012 \\
\end{longtable}

\setlength{\tabcolsep}{4pt}


\begin{longtable}{|c|c|r|r|p{11.4cm}|}
\hline
\textbf{Knot} & \textbf{ID}  &  {$\boldsymbol{\tb}$} & \ {$\boldsymbol{r}$} & \textbf{Parents} \\ \hline
\endhead
            8a1 & S1131 & -8 & -1 & G1011 \\ \cline{2-5}
~ & S1132 & -8 & 1 & G1011 \\ \hline
m8a1 & S2261 & -4 & -1 & G1850, G1849, G1848 \\ \cline{2-5}
~ & S2262 & -4 & 1 & G1849, G1848, G1850 \\ \hline
8a2 & S1232 & -12 & -5 & G1015, G1014 \\ \cline{2-5}
~ & S1231 & -12 & -3 & G1015, G1012, G1014, G1013 \\ \cline{2-5}
~ & S1230 & -12 & -1 & G1017, G1012, G1018, G1013, G1016 \\ \cline{2-5}
~ & S1233 & -12 & 1 & G1018, G1020, G1019, G1017, G1016 \\
\end{longtable}


\begin{longtable}{|c|c|r|r|p{11.4cm}|}
\hline
\textbf{Knot} & \textbf{ID}  &  {$\boldsymbol{\tb}$} & \ {$\boldsymbol{r}$} & \textbf{Parents} \\ \hline
\endhead
8a1 & T1408 & -9 & -2 & S1131 \\ \cline{2-5}
~ & T1409 & -9 & 0 & S1132, S1131 \\ \cline{2-5}
~ & T1410 & -9 & 2 & S1132 \\ \hline
m8a1 & T2846 & -5 & -2 & S2261 \\ \cline{2-5}
~ & T2847 & -5 & 0 & S2261, S2262 \\ \cline{2-5}
~ & T2848 & -5 & 2 & S2262 \\ \hline
8a2 & T1521 & -13 & -6 & S1232 \\ \cline{2-5}
~ & T1520 & -13 & -4 & S1232, S1231 \\
\end{longtable}

\bibliographystyle{mwamsalphack}
\bibliography{references}

\end{document}